\input amstex
\documentstyle{amsppt}
\magnification=\magstep1
\hoffset=.25truein
\vsize=8.75truein
\NoBlackBoxes
\nologo

\define\bs{{\backslash}}
\define\cir{{\, \scriptstyle{\circ}\, }}
\define\IN{{\Bbb N}}
\define\IQ{{\Bbb Q}}
\define\IR{{\Bbb R}}

\define\cC{{\Cal C}}
\define\cI{{\Cal I}}
\define\cJ{{\Cal J}}
\define\cO{{\Cal O}}
\define\cP{{\Cal P}}
\define\hatf{{\widehat f}}
\define\hatg{{\widehat g}}
\define\hatvp{{\widehat\varphi}}
\define\hatxi{{\widehat\xi}}
\define\[[{{[\![}}
\define\]]{{]\!]}}

\define\bcC{{\underline\cC}}
\define\bfm{{\underline m}}

\define\tilC{{\tilde C}}
\define\tilx{{\tilde x}}
\define\tily{{\tilde y}}
\define\tilsig{{\tilde \sigma}}
\define\endprf{{\hfill $\square$}}

\define\codim{{\text{codim}\, }}
\define\exc{{\text{exc}\, }}
\define\grad{{\text{grad}\, }}
\define\Inf{{\text{Inf}\, }} 
\define\inv{{\text{inv}}}
\define\Sing{{\text{Sing}\, }}
\define\supp{{\text{supp}\, }}

\topmatter
\title{Resolution of singularities in Denjoy-Carleman classes}\endtitle
\author{Edward Bierstone and Pierre D. Milman}\endauthor
\rightheadtext{Bierstone and Milman}
\leftheadtext{Resolution of singularities in Denjoy-Carleman classes}
\endtopmatter

\address{E.B. and P.M.: Department of Mathematics, University of
Toronto, Toronto, Ontario, Canada M5S 3G3}\endaddress
\email{bierston\@math.toronto.edu\quad and \
milman\@math.toronto.edu}\endemail

\thanks
The authors' research was partially supported by
NSERC grants OGP0009070, OGP0008949 and the Killam Foundation.
\endthanks

\subjclass{Primary 26E10, 32S45, 58C25; Secondary 14E15,
14P15, 30D60}\endsubjclass
\abstract
We show that a version of the desingularization theorem of Hironaka
holds for certain classes of $\cC^\infty$ functions (essentially, for
subrings that exclude flat functions and are closed under differentiation
and the solution of implicit equations). Examples are quasianalytic classes,
introduced by E. Borel a century ago and characterized by the Denjoy-Carleman
theorem. These classes have been poorly understood in dimension $> 1$.
Resolution of singularities can be used to obtain many new results; for
example, topological Noetherianity, \L ojasiewicz inequalities, division
properties.
\endabstract

\endtopmatter
\document
\baselineskip=18pt 

\head 1. Introduction\endhead

We show that a version of the desingularization theorem 
of Hironaka \cite{Hi1} holds for certain classes of 
infinitely differentiable functions -- essentially, for 
subrings that exclude flat functions and are closed 
under differentiation and the solution of equations
satisfying the conditions of the implicit function
theorem.  Examples are ``quasianalytic classes'', 
introduced by E. Borel a century ago \cite{Bo1} 
and characterized (following questions of Hadamard in studies 
of linear partial equations \cite{Ha}) 
by the Denjoy-Carleman
theorem \cite{De}, \cite{Ca}. (See Section 2 below.)

Quasianalytic classes in one variable play an important part
in harmonic analysis and other areas. 
(See, for instance, \cite{HJ}, \cite{Ko}, \cite{T}.)
In several variables, there are beautiful modern developments
of E.M. Dyn'kin \cite{Dy1}, \cite{Dy2}, but the subject is 
much less understood, perhaps because of a lack of the 
standard techniques of local analytic geometry.  For example, 
the Weierstrass preparation theorem fails 
(Childress \cite{Ch}) and it seems unknown (and unlikely)
that, in general, a ring of germs of functions in a 
Denjoy-Carleman class in Noetherian. 

It may seem surprising that desingularization theorems
nevertheless hold for Denjoy-Carleman classes.  
Our proof of resolution of singularities in \cite{BM2}, 
however, (at least in the case of a hypersurface or 
``principalization of an ideal'' 
\cite{BM2, Thm. 1.10}) uses only elementary 
``differential calculus'' properties that are satisfied 
in these classes.  This was pointed out in 
\cite{BM2, (0.1)} and a simple version of resolution 
of singularities (as in \cite{BM1, Sect. 4}) for quasianalytic 
classes has already been used by Rolin, Speissegger and Wilkie
in their study of $o$-minimality of Denjoy-Carleman
classes \cite{RSW}.

In this article, we isolate the properties of a class of 
$\cC^\infty$ functions needed for resolution of singularities
(Section 3).  We formulate the most general version of 
desingularization known for these classes in Section 5 
(Theorems 5.9 and 5.10; see also Remark 7.10).
Detailed proofs can be found in \cite{BM2}, but we include
a complete proof of a simple version 
(Theorem 5.12; cf. \cite{BM1, Sect. 4}) that in general 
suffices for applications, in a language that makes it 
clear that only the properties of Section 3 are involved. 
The proof (presented in Section 7) is meant at the same 
time to serve as an introduction to two further articles, 
{\it Desingularization algorithms I.  Role of exceptional 
divisors} and {\it Desingularization algorithms II. Binomial 
varieties}, that we plan to publish shortly.

The properties of Section 3 are known for Denjoy-Carleman
classes.  (See references in Section 4.)  We include 
proofs in Section 4 because the minimal hypotheses 
required are not always clear in the literature, nor is the 
elementary nature of the results and the background required 
to prove them. 

Resolution of singularities is of course a powerful tool; 
it can be used to prove several other new results about 
Denjoy-Carleman classes.  Many of the geometric properties
of semialgebraic sets, for example, are satisfied by 
$o$-minimal structures in general \cite{vdDM}.  
The following are discussed in Section 6 below: 
(1) Topological Noetherianity (Theorem 6.1). 
(2) \L ojasiewicz inequalities (Theorem 6.2). 
Proofs of \L ojasiewicz's inequalities depending only on
a simple version of resolution of singularities were 
already given in \cite{BM2, Sect.2}. 
These inequalities were known previously for Denjoy-Carleman 
classes only 
in dimension $2$, under more restrictive hypotheses 
(Vol'berg \cite{V}).  
(3)  Division properties (Theorem 6.3). 

Several unresolved questions about Denjoy-Carleman classes
are raised in the text.  We are grateful to Vincent 
Thilliez for clarifying many points about 
quasianalytic functions.

\head 2. Quasianalytic functions\endhead 

A ``quasianalytic'' function means (roughly speaking) 
a $\cC^\infty$ function that is determined by its Taylor 
expansion at any point.  Quasianalytic functions originate 
in E. Borel's ideas on generalization of the principal of analytic 
continuation.  Borel showed that, if a sequence of complex 
number $\{ A_k\}$ converges to $0$ fast enough, then a 
series $\sum_{k=1}^\infty A_k/(z-a_k)$ converges
normally together with all its derivatives on a ``big''
set of real line segments in a compact set. 
If the poles $a_k$ accumulate everywhere on such a line 
segment, we get a quasianalytic function on the line 
segment that is nowhere analytic \cite{Bo2}.

Let $\cC^\infty (U)$ denote the ring of $\cC^\infty$ 
functions on an open subset $U$ of $\IR^n$. 
Let $f\in \cC^\infty (U)$.  For every $\alpha \in \IN^n$, 
$\alpha = (\alpha_1,\dots, \alpha_n)$, we write 
$$
f_\alpha \ := \ \frac1{\alpha!} D^\alpha f \ , 
$$
where $\alpha! = \alpha_1! \cdots \alpha_n !$ and $D^\alpha$
denotes the partial derivative 
$\partial^{|\alpha|} / \partial x^{\alpha_1}_1 \cdots 
\partial x^{\alpha_n}_n $, 
$|\alpha | = \alpha_1 + \cdots + \alpha_n $. 
($\IN$ denotes the nonnegative integers.)

Let $m = \{ m_0 , m_1 ,\dots \} $ denote a sequence of 
positive numbers. 

\definition{Definition 2.1}
$\cC_m (U) := \{ f \in \cC^\infty (U) : $ for 
every compact $K \subset U$, there are constants 
$A$, $B>0$ such that 
$$
| f_ \alpha (x) | \ \le \ AB^{|\alpha | } m_{|\alpha|} \ , 
$$
for all $\alpha \in \IN^n$ and $x\in K$. 
\enddefinition

Hadamard raised the question of characterizing the sequences
$m$ such that the class $\cC_m$ is {\it quasianalytic}; 
i.e., such that, if $U$ is connected, then the Taylor series
homomorphism
$$
f \ \mapsto \ \hatf_a (x) \ := \ \sum\limits_{\alpha \in \IN^n} 
f_a (a) x^\alpha 
$$
from $\cC_m (U)$ to the ring of formal power series in $n$
indeterminates, in injective for any $a\in U$ 
\cite{Ha, Bk. I, Ch. II}.  The Denjoy-Carleman 
theorem \cite{De}, \cite{Ca} 
is a solution of Hadamard's problem.

We assume that $m = \{ m_k\}$ satifies the hypothesis
$$
\{ m_k \} {\hbox{ is logarithmically convex, }}
\tag 2.2
$$
or, equivalently, 
$$
\left\{ \frac{m_{k+1}}{m_k} \right\} 
{\hbox{ is increasing. }}
$$
(By ``increasing'', we mean ``nondeceasing''; i.e., 
``$\le$''.) The hypothesis (1.1) implies that 
$$
m_j m_k \ \le \ m_o m_{j+k} \ , \quad {\hbox{ for all }} j,k \ , 
$$
and that 
$$
\{ m^{1/k}_k \} {\hbox{ is increasing. }}
$$
The first of these conditions implies that $\cC_m (U)$ is 
a ring, and the second that $\cC_m (U)$ contains the ring 
$\cO (U)$ of real-analytic functions on $U$, for all open 
$U \in \IR^n$. 

Under the hypothesis (2.2), the Denjoy-Carleman theorem 
(see \cite{H\"o, Thm. 1.3.8}, \cite{Ru, Thm. 19.11}) asserts that 
$\cC_m$ is quasianalytic if and only if 
$$
\sum\limits_{k=0}^\infty \ \frac{m_k}{(k+1)m_k} \ = \ \infty \ . 
\tag 2.3
$$
\noindent {\it Note.} In the literature, $\cC_m$ is more usually denoted
$\cC_M$, where $M = \{ M_k \}$ and each $M_k = k! m_k$.  
The use of $\{ m_k\}$ instead of $\{ M_k\}$ and 
$f_\alpha = D^\alpha f / \alpha!$ instead of $D^\alpha f$ will 
be 
convenient for the estimates on derivatives that we make in 
Section 3 below. 

If $m = \{ m_k \}$ satisfies the hypotheses (2.2) and (2.3), 
then $\cC_m$ is called a {\it Denjoy-Carleman class}.

If $f\in \cC_m (U)$, then each partial derivative
$f_{(j)} = \partial f / \partial x_j \in 
\cC_{m^{+1}} (U)$, where $m^{+1}$ denotes the shifted 
sequence
$$
m^{+1} \ := \ \{ m_{k+1} \}_{k\in \IN} \ . 
$$
Clearly, if $m$ satisfies (2.2) (respectively, (2.3)), 
then $m^{+1} $ satisfies (2.2) (respectively, (2.3)). 
If $m = \{ m_k \}$ satisfies the hypothesis
$$
\sup \left( \frac{m_{k+1}}{m_k} \right)^{1/k} \ < \ \infty \ , 
\tag 2.4
$$
then $\cC_{m^{+1}} = \cC_m$, so that $\cC_m$ is closed under 
differentiation.  (Conversely, if $\cC_m$ is closed under 
differentiation, then (2.4) holds (S. Mandelbrojt \cite{M})). 

\head 3. $\cC^\infty$ classes that admit resolution of 
singularities\endhead 

Suppose that, for every open subset $U$ of $\IR^n$, $n\in \IN$, 
we have an $\IR$-subalgebra $\cC (U)$ of $\cC^\infty (U)$. 
Our desingularization thoerems require only the following 
assumptions (3.1)--(3.6) on $\cC(U)$, for any open 
$U \in \IR^n$. 
\smallskip

\noindent (3.1)\qquad $\cP (U)\subset \cC(U)$, where $\cP(U)$
denotes the algebra of restrictions to $U$ of polynomial functions
on $\IR^n$. 
\smallskip

\noindent (3.2)\qquad $\cC$ {\it is closed under composition}.
Suppose that $V$ is an open subset of $\IR^p$ and that 
$\varphi = (\varphi_1 ,\dots, \varphi_p) : U \to V$ is a mapping 
such that each $\varphi_i \in \cC(U)$.  
Then $g\cir \varphi \in  \cC(U)$, for all $g\in \cC(V)$.
\smallskip

A mapping $\varphi : U\to V$ will be called a $\cC$-{\it mapping}
if $g \cir \varphi \in \cC(U)$, for every $g\in \cC(V)$. 
Write $\varphi = (\varphi_1 ,\dots, \varphi_p)$.  
It follows from (3.1) and (3.2) that $\varphi$ is a 
$\cC$-mapping if and only if $\varphi_i \in \cC(U)$, 
$i=1,\dots, p$. 
\smallskip

\noindent (3.3)\qquad $\cC$ {\it is closed under differentiation}.
For all $f \in \cC(U)$, 
$$
\frac{\partial f}{\partial x_i} \ \in \ \cC (U) \ , \quad 
i = 1,\dots, n. 
$$
\smallskip

\noindent (3.4)\qquad $\cC$ {\it is quasianalytic}.
If $f\in \cC (U)$ and $\hatf_a = 0$, where $a \in U$, 
then $f$ vanishes in a neighbourhood of $a$. 
\smallskip

Since $\{ x : \hatf_x = 0 \}$ is closed in $U$, (3.4) is
equivalent to the following property:  If $U$ is connected, then, 
for all $a\in U$, the Taylor series homomorphism 
$\cC (U) \ni f \mapsto \hatf_a \in \IR 
\[[ x\]] $ is injective. 
($\IR \[[ x \]] = \IR \[[ x_1,\dots , x_n \]]$ denotes
the ring of formal power series in $x = (x_1,\dots, x_n)$
with coefficients in $\IR$. 
\smallskip

\noindent (3.5)\qquad $\cC$ {\it is closed under division by 
a coordinate}.  If $f \in \cC (U)$ and 
$$ 
f(x_1,\dots, x_{i-1} , a_i , x_{i+1} ,\dots, x_n ) \ \equiv \ 0 \ , 
$$
then $f(x) = (x_i - a_i ) h(x) $, where $h \in \cC(U)$. 

\smallskip

\noindent (3.6)\qquad $\cC$ {\it is closed under inverse}.
Let $\varphi: U \to V$ be a $\cC$-mapping between open 
subsets $U$, $V$ of $\IR^n$.  Suppose that $a\in U$, 
$\varphi (a) = b$ and the Jacobian matrix 
$$ 
\frac{\partial\varphi }{\partial x} (a) \ := \ 
\frac{\partial (\varphi_1 ,\dots, \varphi_n)}{\partial (x_1,\dots,x_n)}
(a) 
$$
is invertible.  Then there are neighbourhoods 
$U'$ of $a$, $V'$ of $b$, and a $\cC$-mapping 
$\psi: V' \to U'$ such that $\psi (b) = a$ and 
$\varphi \cir \psi$ is the identity mapping of $V'$.
\smallskip

Property (3.6) is equivalent to the following 
{\it implicit function theorem in} $\cC$: 
Suppose that $U$ is open in $\IR^n \times \IR^p$
(with product coordinates 
$(x,y) = (x_1,\dots, x_n $, $ y_1 ,\dots, y_p)$. 
Suppose that $f_1 ,\dots, f_p \in \cC (U)$, 
$(a,b) \in U$, $f(a,b) = 0$ and \break
$(\partial f / \partial y ) (a,b) $ is invertible, where 
$f = (f_1,\dots, f_p)$.  Then there is a product neighbourhood
$V \times W$ of $(a,b)$ in $U$ and a $\cC$-mapping 
$g: V\to W$ such that $g(a) = b$ and 
$$ 
f\big( x,g(x)\big) \ = \ 0 \ ,\quad x \in V . 
$$

Property (3.6) implies that $\cC$ {\it is closed under
reciprocal}; i.e., if $f \in \cC (U)$ vanishes nowhere in 
$U$, then $1/f \in \cC (U)$. 

Under the conditions (3.1)-(3.6), we can use open subsets
$U$ of $\IR^n$ and the algebras of functions $\cC (U)$
as local models in order to define a category $\bcC$
of $\cC$-{\it manifolds} and $\cC$-{\it mappings}. 
The dimension theory of $\bcC$ follows from that of 
$\cC^\infty$ manifolds. We will need two fundamental
properties of such a category $\bcC$: 

\proclaim{Proposition 3.7}
A {\rm smooth} subset of a $\cC$-manifold is a $\cC$-submanifold.
In other words:  Let $M$ be a $\cC$-manifold. 
Suppose that $U$ is open in $M$, $g_1,\dots, g_p \in \cC(U)$, 
and the gradients {\rm grad}$\, g_i$ are linearly independent 
at every point of the zero set 
$$ 
X \ := \ \{ x\in U : \ g_i (x) = 0 \ , \quad 
i = 1,\dots, p\} \ . 
$$
Then $X$ is a closed $\cC$-submanifold of $U$, of 
codimension $p$.
\endproclaim

Proposition 3.7 is of course a consequence of the 
implicit function property (3.6). 

\proclaim{Proposition 3.8}
$\bcC$ is closed under blowing up with centre a closed 
$\cC$-submanifold. 
\endproclaim

\definition{Definition 3.9}
A {\it blowing-up} of a $\cC^\infty$ manifold $M$ with 
{\it centre} a smooth closed subset $C$ is a $\cC^\infty$
mapping $\sigma : M' \to M$ from a $\cC^\infty$ manifold
$M'$, that can be described in local coordinates as 
follows. (See also \cite{BM2, p. 236}.)  First suppose that 
$U$ is an open neighbourhood of $0$ in $\IR^n$, and that $C$ 
is a coordinate subspace
$$ C \ = \ \{ x_i = 0 \ , \quad i \in I \} \ , $$
where $I \subset \{ 1,\dots, n\}$.  The blowing-up 
$\sigma : U' \to U$ with centre $C$ is a mapping where 
$U'$ can be covered by coordinate charts $U_i$, $i\in I$, 
and each $U_i$ has a coordinate system $(y_1,\dots, y_n)$ 
in which $\sigma$ is given by the formulas
$$
\aligned
x_i &= y_i \ , \\
x_j &= y_i y_j \ , \,\quad j\in I \bs \{ i \} \ , \\
x_j &= y_j \ , \qquad j \not\in I \ . 
\endaligned
$$
In general, if $M$ is a $\cC^\infty$ manifold and $C$ 
a closed $\cC^\infty$ submanifold of $M$, then every 
point of $C$ admits a coordinate neighbourhood $U$ in 
which $C$ is a coordinate subspace as above; over this 
neighbourhood, the blowing-up $\sigma : M ' \to M$ 
identifies with the mapping $U' \to U$ defined above, 
On the other hand, $\sigma$ is a diffeomorphism over 
$M\bs C$. The preceding conditions determine
$\sigma : M' \to M$ uniquely, up to a diffeomorphism of 
$M'$ commuting with the projections to $M$. 
\enddefinition

It is easy to see that if $M$ is a $\cC$-manifold and $C$
is a closed $\cC$-submanifold of $M$, then the blowing up 
$\sigma: M' \to M$ with centre $C$ is a $\cC$-mapping. 
This is the assertion of Proposition 3.8.

\head 4. Denjoy-Carleman classes\endhead

Let $m = \{ m_k \}_{k\in \IN}$ denote a sequence of positive 
numbers satisfying (2.2); i.e., $m$ is logarithmically convex. 
Since $\cO (U) \subset \cC_m (U)$ for all open subsets $U$
of $\IR^n$, $n\in \IN$, then $\cC = \cC_m$ satisfies property 
(3.1).  We will show that $\cC = \cC_m$ satisfies properties
(3.2) and (3.6) below (Theorems 4.7 and 4.10).  
The following weaker version of (3.3) is obvious. 

\smallskip

\noindent (3.3$'$)\qquad If $U$ is open in $\IR^n$ and 
$f\in \cC_m (U)$, then each partial derivative 
$\partial f / \partial x_i \in \cC_{m^{+1}} (U)$. 
\smallskip

\noindent By the standard integral formula, we get the 
following weaker version of (3.5). 

\smallskip

\noindent (3.5$'$)\qquad 
If $f\in \cC_m (U)$ and $f(x_1,\dots, x_{i-1}, a_i, x_{i+1},\dots, 
x_n) \equiv 0$, then $f(x) = (x_i - a_i ) h(x) $, where 
$h\in \cC_{m^{+1}} (U)$. 
\smallskip

Therefore, if $\cC = \cC_m$, where $m$ satisfies 
(2.4), or, more generally, if 
$$
\cC \ = \ \bigcup\limits_{j=0}^\infty \ \cC_{m^{+j}} \ , 
\tag 4.1
$$
where $m^{+j}$ denotes the shifted sequence
$$
m^{+j} \ = \ \{ m_{k+j} \}_{k\in \IN} \ , 
$$
then $\cC $ satisfies (3.1), (3.2), (3.3), (3.5) and (3.6). 
Of course, if $m$ satisfies the Denjoy-Carleman condition
(2.3), then $\cC = \cC_m$ and $\cC = \bigcup \cC_{m^{+j}}$
satisfy property (3.4). 

Our proofs of properties (3.2) and (3.6) are based on a 
several-variable version of Fa\`a de Bruno's formula \cite{FdB}.
Consider a composite function $h = f\cir g $, where 
$g(x) = \big(g_1 (x) ,\dots,g_p (x) \big)$, 
$x= (x_1,\dots, x_n)$, and $f(y) = f(y_1,\dots,y_p)$. 
Recall that $f_\alpha (y)$ denotes 
$D^\alpha f(y) / \alpha!$, $\alpha \in \IN^p$.
Write $g_\gamma := (g_{1,\gamma} , \dots, g_{p,\gamma})$, 
$\gamma \in \IN^n$.  Thus the $h_\gamma (x) = (f\cir g)_\gamma
(x)$,  $\gamma \in \IN^n$, are the coefficients of the power 
series in $u$, 
$$ 
\sum\limits_{\gamma \in \IN^n} h_\gamma (x) u^\gamma \ , 
$$
obtained by substituting the power series
$$ 
\sum\limits_{\delta \in \IN^n \bs \{ 0\} } 
g_\delta (x) u^\delta 
$$
for $z = (z_1,\dots, z_p)$ in the power series
$$ 
\sum\limits_{\alpha \in \IN^p } 
f_\alpha \big( g(x) \big) z^\alpha \ . 
$$

\proclaim{Lemma 4.2}
Let $a_i \in \IR^p$, $i=1,\dots, \ell$, and let 
$\alpha = (\alpha_1,\dots, \alpha_p ) \in \IN^p$. Then 
$$ 
(a_1 + a_2 + \dots + a_\ell )^\alpha \ = 
\ \sum \frac{\alpha!}{k_1 ! \cdots k_\ell ! } 
a_1^{k_1} \cdots a^{k_\ell}_\ell \ , 
$$
where the sum is taken over all
$(k_1,\dots, k_\ell)\in (\IN^p)^\ell$ such that 
$\alpha = \sum^\ell_{i=1} k_i $. 
\endproclaim 

\demo{Proof}
$$
\aligned
(a_1 + a_2 + \cdots + a_\ell)^\alpha 
\ &= \ \prod_{j=1}^p (a_{1j} + \cdots + a_{\ell j} )^{\alpha_j}\\
\ &= \ \prod_{j=1}^p \left( 
\sum_{\alpha_j = \sum k_{ij},\atop k_{ij} \in \IN } 
\frac{\alpha_j !}{k_{1j}! \cdots k_{\ell j}! } 
a^{k_{1j}}_{1j} \cdots a^{k_{\ell j}}_{\ell j} 
\right) \ . 
\endaligned
$$
In the expansion of the latter product, each term is a unique 
product of terms, one from each of the $p$ factors in the 
product. 
\endprf
\enddemo

\proclaim{Proposition 4.3}
(Fa\`a de Bruno's formula in several variables.)
For all $\gamma \in \IN^n \bs \{0\}$, 
$$
h_\gamma (x) \ = \ \sum 
\frac{\alpha!}{k_1 ! \cdots k_\ell !} 
f_\alpha \big( g(x)\big) 
g_{\delta_1} (x)^{k_1} \cdots g_{\delta_\ell} (x)^{k_\ell} \ , 
$$
where $\alpha = k_1 + \cdots + k_\ell$ and the sum is taken 
over all sets $\{ \delta_1,\dots , \delta_\ell \}$ of 
$\ell$ distinct elements of $\IN^n \bs \{ 0 \}$
and all ordered $\ell$-tuples $(k_1,\dots , k_\ell) \in 
(\IN^p \bs \{0 \} )^\ell $, $\ell = 1, 2, 3, \dots , $ such that 
$$
\gamma \ = \ \sum_{i=1}^\ell | k_i | \delta_i \ . 
$$
\endproclaim 

\demo{Proof}
$$
\sum_{\gamma\in \IN^n} h_\gamma (x) u^\gamma \ = 
\ \sum_{\alpha\in \IN^p} f_\alpha \big(g(x)\big)
\left( \sum_{\delta \in \IN^n \bs \{0 \}}
g_\delta (x) u^\delta \right)^\alpha \ , 
$$
so the result follows immediately from Lemma 4.2. 
\endprf
\enddemo

In the remainder of this section, $m = \{ m_k \}_{k\in \IN}$
denotes a logarithmically convex sequence of positive numbers. 
The following inequality of Childress \cite{Ch} 
is obviously connected to the Fa\`a de Bruno formula in 
one variable. 

\proclaim{Proposition 4.4}
Let $k_1,\dots, k_n$ be nonnegative integers such that 
$k_1 + 2k_2 + \cdots + nk_n = n$.  Set $k=k_1+\cdots + k_n$. 
Then 
$$
m_k m^{k_1}_1 \cdots m^{k_n}_n \ \le \ m^k_1 m_n \ . 
$$
\endproclaim 

\demo{Proof}
The result is trivial if $k_n = 1$; we can therefore assume 
that $k_n = 0$. 

{\it Case I.} \ $k_1 \ne 0$. Let $k'_1 = k_1 - 1$, 
$k'=k-1$. Then $k' = k'_1 + k_2 + \cdots + k_{n-1} $ and 
$n-1 = k'_1 + 2k_2 + \cdots + (n-1) k_{n-1}$. By induction on $n$, 
$$
m_{k-1} 
m^{k'_1}_1 m^{k_2}_2  \cdots 
m^{k_{n-1}}_{n-1} m^{k_n}_n \ \le \ m^{k'}_1 m_{n-1}
$$
(remember $m^{k_n}_n = 1 ! $);  thus 
$$
\aligned
m_k m^{k_1}_1 \cdots m^{k_n}_n 
\ &= \ m_1\frac{m_k}{m_{k-1}} m_{k-1} m^{k'_1}_1 m^{k_2}_2 
\cdots m^{k_n}_n \\
&= \ m_1\frac{m_n}{m_{n-1}} m^{k'}_1 m_{n-1} \\
&= \ m^k_1 m_n \ . 
\endaligned
$$

{\it Case II.} \ $k_1 = 0$.  We have 
$$ 
n-k \ = \ k_2 + 2k_3 + \cdots + (n-k) k_{n-k+1} + \cdots ; 
$$
thus $k_j = 0$ if $j > n-k+1$, and $k= k_2 + \cdots + k_{n-k+1}$. 
By induction, 
$$
m_{k+1} m^{k_2}_2 \cdots m^{k_{n-k+1}}_{n-k+1} \ \le 
\ m^k_2 m_{n-k+1} \ ; 
$$
in other words, 
$$
m_{k+1} m^{k_1}_1 \cdots m^{k_n}_n \ \le \ m^k_2 m_{n-k+1} \ . 
$$
Therefore, 
$$
\aligned
m_k m^{k_1}_1 \cdots m^{k_n}_n 
\ &\le \ \frac{m_k}{m_{k+1}} m^k_2 m_{n-k+1}\\
&\le \ m_1 m^{k-1}_2 m_{n-k+1}\\
&\le \ m^2_1 m^{k-2}_2 m_{n-k+2}\\
&\le \ \cdots \ \le \ m^k_1 m_n \ .
\endaligned
$$
{\vskip -.30truein}\line{\hfill \endprf}
\enddemo

\proclaim{Corollary 4.5}
Let $k_1 ,\dots k_\ell \in \IN^p \bs \{0 \}$ and 
$\delta_1 ,\dots , \delta_\ell \in \IN^n \bs \{0 \}$. 
Set $\alpha = k_1 + \cdots + k_\ell$ and 
$\gamma = |k_1| \delta_1 + \cdots + |k_\ell | \delta_\ell $. 
Then 
$$
m_{|\alpha| } 
m^{|k_1|}_{|\delta_1|} \cdots 
m^{|k_\ell|}_{|\delta_\ell |} \ \le 
\ m^{|\alpha|}_1 m_{|\gamma|} \ . 
\tag 4.6
$$
\endproclaim 

\demo{Proof}
This is a special case of Childress's inequality because the 
latter applies with some $k_i = 0$. 
(We can assume that all $|\delta_i|$ are distinct because if
$|\delta_i | = |\delta_j|$ for some $i$ and $j\ne i$, then 
we can replace  
$m^{|k_i|}_{|\delta_i|}  m^{|k_j|}_{|\delta_j|} $ in the 
left-hand side of (4.6) by 
$m^{|k_i| + |k_j|}_{|\delta_i|}$.) 
\endprf
\enddemo

\proclaim{Theorem 4.7} (Composition; cf. \cite{Rou}.)
Let $U$ and $V$ denote open subsets of $\IR^n$ and $\IR^p$, 
respectively.  Let $f\in \cC_m (V)$ and let $g = (g_1,\dots, 
g_p) : U \to V$, where each $g_j \in \cC_m (U)$.  
Then $f\cir g \in \cC_m (U)$. 
\endproclaim 

\demo{Proof}
Let $K$ be a compact subset of $U$.  Then there are 
constants $a$, $b$, $c$, $d > 0$ such that 
$$
\aligned
|f_\alpha (y)| \ &\le \ ab^{|\alpha|} m_{|\alpha|} \ , \quad 
\text{for all } y\in g(K) , \alpha \in \IN^p, \\
|g_{j,\delta} (x)| \ &\le \ cd^{|\delta|} m_{|\delta |} \ , \quad 
\ \text{for all } x\in K, \delta \in \IN^p, j=1,\dots , p. 
\endaligned
$$
Let $h= f\cir g$. Let $\gamma \in \IN^n \bs \{ 0\}$. 
By Proposition 4.3 and Corollary 4.5, if $x\in K$, then 
$$
\aligned 
|h_\gamma (x)| 
\ &\le \ a \sum \frac{\alpha !}{k_1 ! \cdots k_\ell !} 
(bc)^{|\alpha | } \delta^{| \gamma |} 
m_{|\alpha|} m^{|k_1|}_{|\delta_1 |} \cdots 
m^{|k_\ell |}_{|\delta_\ell |} \\
&\le \ ad^{|\gamma|} m_{|\gamma|} 
\sum \frac{\alpha !}{k_1 ! \cdots k_\ell !} 
(bc m_1 )^{|\alpha|} \ . 
\endaligned 
$$
By Lemma 4.8 following, there are constants
$C$, $D$ depending only on $bc m_1$, $n$ and $p$, such that 
$$
\sum \frac{\alpha !}{k_1 ! \cdots k_\ell !} 
(bc m_1)^{|\alpha|} \ \le \ CD^{|\gamma |} \ , 
$$
for all $\gamma \in \IN^n \bs \{ 0\} $. 
(The summation is always as in Proposition 4.3.)
Thus, 
$$
|h_\gamma (x) | \ \le \ aC (dD)^{|\gamma|} m_{|\gamma |} \ . 
$$
{\vskip -.30truein}\line{\hfill \endprf}
\enddemo

\proclaim{Lemma 4.8}
Let $\lambda > 0$.  Set $H(u) = \sum_{\gamma\in \IN^n}
H_\gamma u^\gamma $, $u= (u_1,\dots, u_n)$, 
where $H_0 = 1$ and, for each $\gamma \in \IN^n \bs \{0\} $, 
$$ 
H_\gamma = \sum \frac{\alpha!}{k_1 ! \cdots k_\ell !} 
\lambda^{|\alpha |} 
$$
(summation as in Proposition 4.3).  Then $H$ is a convergent
power series. 
\endproclaim 

\demo{Proof}
Define 
$$
\aligned
G_j (u_1,\dots, u_n ) \ &:= \ \prod^n_{i=1} 
\left( \frac{1}{1-u_i} \right) - 1 \ , \quad j=1,\dots,p , \\
F (z_1,\dots, z_p ) \ &:= \ \prod^p_{j=1} 
\left( \frac{1}{1-\lambda z_j } \right) \ , 
\endaligned
$$
so that 
$$
\aligned
G_j (u) \ &= \ \prod^n_{i=1} (1+u_i + u^2_i + \cdots ) - 1\\
	&= \ \sum_{\delta \in \IN^n \bs \{0 \} } u^\delta \ , \\
F(z) \ &= \ \sum_{\alpha \in \IN^p} (\lambda z)^\alpha \ . 
\endaligned
$$
Then $H = F\cir G$, by Proposition 4.3. 
\endprf
\enddemo

\remark{Remark 4.9}
In the $1$-variable case of Theorem 4.7 $(n=p=1)$, 
we can use Fa\`a de Bruno's formula to show that the constants
$C$ and $D$ in the proof can be taken more precisely as 
$C = bc m_1 $, $D= 1+bcm_1 $
(cf, \cite{KP, Proposition 1.3.3}).
\endremark

\proclaim{Theorem 4.10} (Inverse function theorem; cf.
\cite{Kom}.)
Let $f:U\to V$ denote a $\cC_m$-mapping between open 
subsets $U$, $V$ of $\IR^n$.  Let $x_0 \in U$.  Suppose 
that the Jacobian matrix 
$(\partial \varphi / \partial x) (x_0)$ is invertible.
Then there are neighbourhoods $U'$ of $x_0$, $V'$ of 
$y_0 := f(x_0)$ and a $\cC_m$-mapping $g: V' \to U'$ such that 
$g(y_0) = x_0 $ and $f\cir g$ is the identity mapping 
of $V'$. 
\endproclaim

\demo{Proof}
Write $f = (f_1,\dots, f_n)$.  We can assume that $f$ has 
a $\cC^\infty$ inverse.  Let $K$ be a compact subset of $U$. 
Then there are constants $a$, $b > 0$ such that, for all 
$\alpha \in \IN^n$, $x\in K$ and $i=1,\dots, n$, 
$$
| f_{i,\alpha} (x) | \ \le \ ab^{|\alpha| } m_{|\alpha |} \ . 
$$
Let $x_0 \in K$.  Consider the solution $x = g(y)$ of the equation 
$$ 
f(x_0 + x) \ = \ f(x_0) + y \ . 
$$
We want to show there are constants $c$, $d > 0$ 
{\it independent of} $x_0 \in K$, such that 
$$
|g_{j, \beta} (0) | \ \le \ cd^{|\beta|} m_{|\beta|} \ , 
$$
for all $\beta \in \IN^n $ and $j = 1,\dots, n$. 

Write 
$$
f(x_0 + x) - f(x_0) \ = 
\ \frac{\partial f}{\partial x} (x_0) \cdot x - 
\frac{\partial f}{\partial x} (x_0) \cdot \varphi(x)  \ ; 
$$
in other words, 
$$
\varphi(x) \ = \ x - \frac{\partial f}{\partial x} (x_0)^{-1} 
\big( f(x_0 + x) - f(x_0) \big) \ . 
$$
Set $\Theta (x) := (\partial f / \partial x) (x)^{-1} 
= \big(\theta_{ij} (x)\big) $.  Thus, 
$$
g(y) \ = \ \Theta (x_0) \cdot y + \varphi \big( g(y)\big) \ , 
$$
where $\varphi (0) = 0$, 
$(\partial \varphi / \partial x_k) (0) = 0$, 
$k= 1,\dots, n$, and 
$$
\varphi_{i,\alpha } (0) \ = \ -\sum_{j=1}^n 
\theta_{ij} (x_0 ) f_{j,\alpha } (x_0 ) \ , 
$$
for all $i=1,\dots , n$ and $|a| \ge 2$.  Choose $r>0$
such that 
$$ 
| \theta_{ij} (x) | \ \le \ r 
$$
for all $i$, $j = 1,\dots , n $ and $x \in K$.  Then, 
for all $i$, $\alpha$, 
$$
|\varphi_{i,\alpha} (0) | \ \le 
\ nr ab^{|\alpha|} m_{|\alpha |} \ . 
$$
\enddemo

By Proposition 4.3, if $|\gamma | \ge 2$, then 
$$ g_\gamma (0) \ = \ \sum \frac{\alpha!}{k_1! \cdots k_\ell !}
\varphi_\alpha (0) g_{\delta_1} (0)^{k_1} \cdots g_{\delta_\ell} 
(0)^{k_\ell} \ ,
$$
where only terms with $|\alpha | \ge 2$ are nonzero, 
so only $g_{\delta_j} (0)$ with $|\delta_j | < |\gamma| $ occur.  
(The latter remark is also clear from the method of undetermined
coefficients applied to (4.12).) 

Let $\Phi (x) = \sum \Phi_\alpha x^\alpha$ denote the 
{\it convergent} power series
$$
\Phi (x) \ = \ \sum\limits_{|\alpha| \ge 2} n r a (m_1 b)^{|\alpha|} 
x^\alpha \ , 
$$
and consider the system of equations 
$$
G_i (y) \ = \ \frac{r}{m_1} (y_1 + \cdots + y_n) + 
\Phi (G(y)) \ , \qquad 
i= 1,\dots, n \ ; 
\tag 4.13
$$
write $G_i (y) =\sum G_{i,\gamma} y^\gamma $ and
$G = (G_1,\dots, G_n) $.  Then necessarily $G_1 = \cdots = G_n$;  
the solution of (4.13) is convergent, so there are constants 
$c$, $d$ depending only on $n$, $m_1$, $r$, $a$ and $b$, 
such that 
$$ |G_{i,\gamma} | \ \le \ cd^{|\gamma|} $$
for all $i$, $\gamma$.  Note that {\it all $G_{i,\gamma}$ are 
nonnegative} (recursively by the Fa\`a de Bruno formula, 
or by (4.13)). We claim that 
$$
|g_{i,\gamma} (0) | \ \le \ G_{i,\gamma} m_{|\gamma|} \ ,\qquad 
|\gamma | \ge 1 \ ; 
\tag 4.14
$$
this gives (4.11). 

We prove (4.14) by induction on $|\gamma |$.  To begin, 
consider $\gamma = (j)$ (where $(j)$ denotes the 
multiindex with $1$ in the $j$'th place and $0$ elsewhere): 
$$
|g_{i,(j)} (0)| \ = 
\ |\theta_{ij} (x_0) | \ \le \ r \ = 
\ G_{i,(j)} m_1 \ . 
$$
By Proposition 4.3, Corollary 4.5 and the induction hypothesis, 
if $|\gamma | \ge 2$, then 
$$ \align
|g_{i,\gamma} (0) |
\ &\le \ \sum \frac{\alpha!}{k_1 ! \cdots k_\ell !} 
|\varphi_{i,\alpha} (0) | 
| g_{\delta_1} (0)^{k_1} | \cdots |g_{\delta_\ell} (0)^{k_\ell} | \\
&\le \ \sum \frac{\alpha!}{k_1 ! \cdots k_\ell !} 
n r a b^{|\alpha|} m_{|\alpha|} 
G^{k_1}_{\delta_1} \cdots G^{k_\ell}_{\delta_\ell} 
m^{|k_1|}_{|\delta_1|}  \cdots m^{|k_\ell |}_{|\delta_\ell |} \\
&= \ \sum \frac{\alpha!}{k_1 ! \cdots k_\ell !} 
\Phi_\alpha G^{k_1}_{\delta_1} \cdots G^{k_\ell}_{\delta_\ell} 
\frac{m_{|\alpha|}}{m^{|\alpha|}_1} m^{|k_1|}_{|\delta_1 |}
\cdots m^{|k_\ell|}_{|\delta_\ell|} \\
&\le \ m_{|\gamma|} \sum \frac{\alpha!}{k_1 ! \cdots k_\ell !} 
\Phi_\alpha G^{k_1}_{\delta_1} \cdots G^{k_\ell}_{\delta_\ell} \\
&= \ G_{i,\gamma} m_{|\gamma|} \ , 
\endalign
$$
where the last equality is the Fa\`a de Bruno formula for 
the coefficients of $G_i$, from (4.13). 
\endprf

\remark{Remark 4.15}
We can again get a more precise estimate in the $1$-variable 
case. If $f\in \cC_m$ and $f(g(x))= x$, choose $a$, $b > 0$
such that 
$$\align
|f_{(1)}| \ &\ge \ \frac1{am_1}\\
|f_{(n+1)} | \ &\ge \ \frac1{(n+1) am^2_1 } b^n m_{n+1} \ , \qquad 
n\ge 1
\endalign
$$
(on a compact subset of $\IR$).  Then 
$$ 
|g_{(n)} | \ \le \ (-1)^{(n-1)} 
{1/2\choose n} 2^n a^n c^{n-1} m_n \ , 
$$
where $c = 2bm_1 $ (cf. \cite{KP, Theorem 1.4.3}). 
\endremark

\head 5. Resolution of singularities\endhead

In this section, $\cC$ denotes a class of $\cC^\infty$ functions
satisfying the hypotheses (3.1)-(3.6). 

\subhead Spaces of class $\cC$ \endsubhead
Let $M$ be a manifold of class $\cC$.  Let $\cO^\cC = \cO^\cC_M$
denote the sheaf of germs of functions of class $\cC$ at points
of $M$.  We regard $M$ as a local-ringed space 
$M = (|M|, \cO^\cC_M )$, where $|M|$ denotes the underlying
topological space of $M$.  (Each stalk $\cO^\cC_a$ of 
$\cO^\cC$ is a local ring.) 
We usually do not distinguish in notation between $M$ and 
$|M|$.  If $\dim_a M = n$, then the completion of 
$\cO^\cC_a$ in the Krull topology 
can be identified with the ring of formal 
power series in $n$ indeterminates. 

Let $\cI \subset \cO^\cC$ denote a sheaf of ideals. 
We say that $\cI$ is of {\it finite type} if, for each 
$a\in M$, there is an open neighbourhood $U$ of $a$ and 
finitely many sections $f_1,\dots, f_q \in \cO^\cC (U) = \cC(U)$
such that, for all $b\in U$, the stalk $\cI_b$ is 
generated by the germs $f_{i,b}$ of the $f_i $ at $b$. 

Suppose that $\cI$ is an ideal of finite type in $\cO^\cC$
(i.e., a subsheaf of ideals of finite type).  Let 
$$
|X| \ := \ \supp \frac{\cO^\cC}{\cI} \ , \qquad 
\cO^\cC_X \ := \ \left( \frac{\cO^\cC}{\cI} \right) 
\Big| |X| \ . 
$$
We call $X = (|X| , \cO^\cC_X)$ a ({\it closed})
$\cC$-{\it subspace} of $M$, and $|X|$ a ({\it closed}) 
$\cC$-{\it subset}. 
(We again usually do not distinguish in notation between $X$
and $|X|$.)  Write $\cI = \cI_X$. 

A closed $\cC$-subspace $X$ of $M$ is a {\it hypersurface}
if $\cI_X$ is a principal ideal (i.e., a sheaf of principal 
ideals). 

We say that a closed $\cC$-subspace $X$ of $M$ is {\it smooth at} 
at point $a\in X$ (or that $a$ is a {\it smooth point} of $X$)
if $\cI_{X,a}$ is generated by elements $f_1,\dots, f_q$ whose 
gradients are linearly independent at $a$. 
Let $\Sing X \subset |X|$ denote the complement of the 
set of smooth points.  By Proposition 3.7, a smooth 
$\cC$-subspace of $M$ is a $\cC$-submanifold. 

Suppose that $\varphi : N \to M$ is a $\cC$-mapping 
of $\cC$-manifolds.  If $\cI \subset \cO^\cC_M$ is an ideal of 
finite type, let $\varphi^{-1} (\cI) \subset \cO^\cC_N$ denote 
the ideal sheaf $\varphi^* (\cI) \cdot \cO^\cC_N$ whose stalk
at each point $b\in N$ is generated by the ring of pull-backs
$\varphi^* (\cI)_b$ of all elements of $\cI_{\varphi(b)}$. 
If $X$ is a closed $\cC$-subspace of $M$, let 
$\varphi^{-1} (X)$ denote the closed $\cC$-subspace of $N$
determined by the ideal $\varphi^{-1} (\cI_X)$. 

\subhead Transformations by blowing up \endsubhead
Let $M$ denote a $\cC$-manifold, and $C$ a closed $\cC$-submanifold
of $M$.  Let $\sigma : M' \to M$ denote the blowing-up of $M$ 
with centre $C$ (Definition 3.9).  Then $\sigma^{-1} (C)$ is a smooth 
closed hypersurface in $M'$; we write $y_\exc $ to denote a
generator of $\cI_{\sigma^{-1} (C) , a' }$, at any point 
$a'\in M'$.  Let $\cI \subset \cO^\cC_M$ be a sheaf of ideals of
finite type.

\definition{Definition 5.1}
If $a\in M$, then the {\it order of $\cI$ at $a$}, 
$$
\mu_a (\cI ) \ := \ \max \{ k \in \IN : \ \cI_a \subseteq 
\bfm^k_a \} \ , 
$$
where $\bfm_a = \bfm_{M,a} $ denotes the maximal ideal of
$\cO^\cC_{M,a}$.  If $a\in C$, then the {\it order of $\cI$
along $C$ at $a$, }
$$
\mu_{C,a} (\cI) \ := \ \max \{ k \in \IN : \ \cI_a \subset
\cI^k_{C,a} \} \ . 
$$
(If $g$ is a germ of a function of class $\cC$ at $a$, 
we define $\mu_a (g)$ and $\mu_{C,a} (g)$ in the same way; 
$\mu_a (g) = \mu_a (\cI)$ and $\mu_{C,a} (g) = \mu_{C,a} (\cI)$, 
where $\cI$ is the ideal generated by $g$.)
\enddefinition

\proclaim{Lemma 5.2}
Each point of $M$ admits a neighbourhood $U$ in which 
$\mu_x (\cI)$ takes only finitely many values and, for any 
$d\in\IN$, $Z_d := \{ x\in U : \mu_x (\cI) \ge d \}$ is a 
closed $\cC$-subset of $U$. 
\endproclaim

(We say that $\mu_x (\cI)$ is {\it Zariski-semicontinuous
(relative to the class $\cC$)}; cf. \cite{BM2, Lemma 3.10}.)

\demo{Proof}
Let $a\in M$.  Let $U$ be an open neighbourhood of $a$ for
which there are $g_1 ,\dots, g_q \in \cI(U)$ that 
generate $\cI_x$, for all $x\in U$.  Then, for any $d\in \IN$, 
$Z_d = \{ x\in U: D^\alpha g_i (x) = 0$, 
$|\alpha | < d $, $i= 1,\dots, q\}$. 
After shrinking $U$ if necessary, $\mu_x (\cI) \le \mu_a (\cI)$, 
for all $x\in U$. 
\endprf
\enddemo

\proclaim{Lemma 5.3}
If $a\in C$, then 
$$ 
\mu_{C,a} (\cI) \ = \ \min \{ \mu_x (\cI) : \ x\in C {\hbox{ near }}
a \} \ . 
$$
In particular, $\mu_{C,x} (\cI)$ is locally constant on $C$. 
\endproclaim 

This is a simple exercise. 

\definition{Definitions and Remarks 5.4}
Let $\cI \subset \cO^\cC_M$ be a sheaf of ideals of finite type. 
The {\it (weak) transform $\cI'$ of $\cI$ by} the blowing-up 
$\sigma$ is a sheaf of ideals of finite type in $\cO^\cC_{M'}$
defined as follows: 
Let $a' \in M'$ and $a= \sigma (a')$.  Then 
$\cI'_{a'}$ is the ideal generated by 
$$
g' \ := \ y^{-d}_\exc g \cir \sigma \ , \quad 
g \in \cI_{X,a} \ , 
\tag 5.5
$$
where $d$ is the largest power of $y_\exc$ that factors
from all $g\in \cI_{X, a}$.  
(If $a' \notin \sigma^{-1} (C)$, then we can take 
$y_\exc = 1$ at $a'$, and $g'= g\cir \sigma $ in (5.5).)

It is easy to see that if $a' \in \sigma^{-1} (C)$, then 
$d= \mu_{C,a} (\cI)$;  it follows from Lemma 5.3 that 
$\cI'$ is of finite type. 

Let $X\subset M$ denote a closed $\cC$-subspace and 
let $X'\subset M'$ denote the closed $\cC$-subspace of $M'$
given by $\cI_{X'} := \cI'$, where $\cI'$ is the above 
transform of $\cI := \cI_X$; $X'$ is called the 
{\it weak transform of $X$ by $\sigma$}.

{\it Suppose that $X$ is a hypersurface.}  
In this case, $X'$ is also called the {\it strict transform of 
$X$ by $\sigma$.}  If $a\in C$ and $g$ denotes a generator 
of $\cI_{X,a}$, then, for all $a'\in \sigma^{-1} (a)$, 
$\cI_{X',a'}$ is the ideal generated by 
$$
g' \ := \ y^{-d}_\exc g\cir \sigma \ , 
$$
where $d$ is the largest power of $y_\exc $ that factors 
from $g\cir \sigma$. 

If $\codim C = 1$, then the blowing-up $\sigma$ is the identity 
mapping, but the transforms above still make sense.  For example, 
if $\cI = \cI_C$ (or if $X = C$) , then $\cI' = \cO^\cC_{M'}$
(or $X' = \emptyset$). 
\enddefinition

\remark{Remark 5.6}
In general, the notions of weak and strict transform do not 
coincide;  see \cite{BM2, Section 3}.
For Denjoy-Carleman classes (or spaces of class $\cC$), it is not 
clear that the strict transform $X'$ is, in general, 
even well-defined as a closed $\cC$-subspace of $M'$. 
One proves that the strict transform is well-defined in a 
category of schemes or analytic spaces, for example, using 
Noetherianity of the local rings and ``Oka-Cartan theory'' 
of coherent sheaves.  (See \cite{BM2, Prop. 3.13 ff.}.)
\endremark

\subhead Desingularization theorems \endsubhead
Le $M$ denote a $\cC$-manifold.  Let $\cI \subset \cO^\cC_M$ 
denote a sheaf of ideals of finite type.  We consider 
sequences of transformations 
$$
\matrix 
\longrightarrow &M_{j+1} 
&{\buildrel {\sigma_{j+1}}\over \longrightarrow}
&M_j &\longrightarrow &\cdots &\longrightarrow 
&M_1 &{\buildrel {\sigma_1}\over \longrightarrow} &M_0 &= &W\\
& & & & & & & & & & &\\
&\cI_{j+1} & &\cI_j & & & &\cI_1 & &\cI_0 &= &\cI | W\\
& & & & & & & & & & &\\
&E_{j+1} & &E_j & & & &E_1 & &E_0 &= &\emptyset
\endmatrix 
\tag 5.7
$$
where $W$ is an open subset of $M$ and, for each $j$: 

\noindent (1) $\sigma_{j+1} : M_{j+1} \to M_j$ is a blowing-up 
with smooth centre $C_{j+1} \subset M_j$, $\cI_{j+1} $
is the transform of $\cI_j$ by $\sigma_{j+1}$, 
and $E_{j+1}$ is the collection of 
{\it exceptional hypersurfaces}
$$
E_{j+1} \ := \ E'_j \cup \{ \sigma^{-1}_{j+1} (C_{j+1}) \} \ , 
$$
where $E'_j$ denotes the collection of strict transforms
$H'$ by $\sigma_{j+1}$ of all hypersurfaces $H\in E_j$. 

\noindent (2) $C_{j+1}$ and $E_j$ {\it simultaneously have only 
normal crossings} (i.e., locally, we can choose coordinates 
with respect to which $C_{j+1}$ is a coordinate subspace and 
$E_j$ is a collection of coordinate hyperplanes. 

We say that the blowing-up $\sigma_{j+1}$ (or the 
centres $C_{j+1}$) in (5.7) are {\it admissible} 
(or $\mu$-{\it admissible}) if, in addition, 
$\mu_a (\cI_j)$ is locally constant on $C_{j+1}$, 
for each $j$. 

If $X$ is a closed $\cC$-hypersurface in $M$ and $\cI=\cI_X$, 
then each $\cI_{j+1} $ in (5.7) is the ideal sheaf 
$\cI_{X_{j+1}}$ of the strict transform $X_{j+1}$ of $X_j$ 
(where $X_0 = X | W$).  In this case, we also write 
$\mu_{X,a} := \mu_a (\cI_X)$; $\mu_{X,a}$ is called the 
{\it order of $X$ at $a$}.

The condition (2) in (5.7) guarantees inductively that each 
$E_{j+1}$ has only normal crossings (i.e., locally, we can choose 
coordinates with respect to which every element of $E_{j+1}$ is 
a coordinate subspace), according to the following simple lemma.

\proclaim{Lemma 5.8} 
Let $H_1,\dots, H_q$ denote smooth $\cC$-hypersurfaces in $M$
that simultaneously have only normal crossings.  Let 
$\sigma:M' \to M$ denote a blowing-up with centre $C$ a smooth 
$\cC$-subspace, such that $C$, $H_1,\dots, H_q$ simultaneously 
have only normal crossings.  Then the strict transforms 
$H'_1 ,\dots, H'_q$ together with $H'_{q+1} := \sigma^{-1} (C)$
simultaneously have only normal crossings.
\endproclaim

\proclaim{Theorem 5.9}
Let $M$ denote a $\cC$-manifold and let $\cI \subset \cO^\cC_M$
be a sheaf of ideals of finite type. 
Let $K$ be a compact subset of $M$. Then there is a neighbourhood
$W$ of $K$ and a finite sequence (5.7) of admissible blowings-up
$\sigma_j$, $j=1,\dots, k+1$, such that $\cI_{k+1} = 
\cO^\cC_{M_{k+1}}$ and $\sigma^{-1} (\cI)$ is a normal-crossings
divisor, where $\sigma :M_{k+1} \to W \subset M$ denotes
the composite of the $\sigma_j$. 

In fact, there is a finite sequence (5.7) satisfying the preceding
assertions and the additional condition that, if $\cJ_\sigma$ 
denotes the ideal generated by the Jacobian determinant of 
$\sigma$ (with respect to any local coordinate systems), 
then $\cJ_\sigma \cdot \sigma^{-1} (\cI)$ is a normal-crossings
divisor.
\endproclaim

({\it Normal-crossings divisor} means a principal ideal of 
finite type, generated locally by a monomial in suitable 
coordinates.)

Suppose that $X$ is a closed $\cC$-hypersurface in $M$.
Clearly then, $\Sing X$ is a closed $\cC$-subset of $M$. 
(It is defined locally by a generator $g$ of $\cI_X$ 
together with all first-order partial derivatives of $g$.)

\proclaim{Theorem 5.10}
Let $M$ denote a $\cC$-manifold, and let $X$ be a closed
$\cC$ hypersurface in $M$.  Set $\cI = \cI_X$.  Let $K$
be a compact subset of $M$.  Then there is an open neighbourhood
$W$ of $K$ and a finite sequence (5.7) of admissible blowings-up
$\sigma_j$, $j=1,\dots, k$, such that: \newline
\noindent (1) for each $j=0,\dots, k-1$, either 
$C_{j+1} \subset \Sing X_j$ or $X_j$ is smooth and 
$C_{j+1} \subset X_j \cap E_j $; \newline
\noindent (2) $X_k$ is smooth; \newline
\noindent (3) $X_k$, $E_k$ and the Jacobian ideal $\cJ_\sigma $
(cf. Theorem 5.9) simultaneously have only normal crossings. 
\endproclaim

Theorem 5.9 in the case of a sheaf of principal ideals is an 
immediate consequence of Theorem 5.10. 
(Let $X$ be the $\cC$-hypersurface determined by the ideal 
$\cI$.  Then $X_k$ (from Theorem 5.10) is smooth and of 
codimension $1$.  Let $\sigma_{k+1} $ be the blowing-up 
with centre $C_{k+1} = X_k$.  Then $\sigma_{k+1}$ is the 
identity, but the strict transform $X_{k+1}$ of $X_k$ is 
empty; i.e., $\cI_{k+1} = \cO^\cC_{M_{k+1}}$.) 
Theorems 5.9 and 5.10 are, in fact, proved in \cite{BM2}
using the same {\it desingularization algorithm}; we refer to 
\cite{BM2} for details, but the idea is very roughly as follows:
There is an invariant 
$$
\inv_\cI (a) \ = \ \big( \nu_1 (a) , s_1 (a) ; \dots ; 
\nu_{t+1} (a) \big) \ , \qquad 
a\in M_j \ , 
\tag 5.11
$$
defined recursively over a sequence of transformations (5.7) 
whose successive centres are ``$\inv_\cI$-admissible'', 
where $\nu_1 (a) = \mu_a (\cI_j)$ and $t\le \dim_a M_j $. 
Sequences of the form (5.11) can be ordered lexicographically; 
$\inv_cI (\cdot)$ takes only finitely many values in a 
neighbourhood $W_j$ of the compact subset 
$K_j = \sigma^{-1}_j (K_{j-1})$ (where $K_0 = K$), 
and the maximum locus of $\inv_cI $ is a union of smooth 
closed $\cC$-subsets of $W_j$ having only normal crossings. 
The desingularization algorithm is given by taking as each 
successive centre $C_{j+1}$ one of the components of the maximum 
locus of $\inv_\cI$  on $W_j$.  Theorems 5.9 and 5.10 follow 
from the basic properties of $\inv_\cI$ (given in 
\cite{BM2, Theorem 1.14}). 

The desingularization theorems are proved in \cite{BM2} in a 
language common to algebraic schemes and analytic spaces
(in characteristic zero), or hypersurfaces of class $\cC$ as
here.  (See \cite{BM2, (0.1)}.)  
In Section 7, we present a complete proof of a simple version 
of Theorem 5.10 (Theorem 5.12 below) that suffices for applications
of the kind considered in Section 6 (or, for example, 
in \cite{RSW}). 
The proof is similar to that of \cite{BM1, Theorem 4.4}, 
which is the source of the desingularization algorithms in 
\cite{BM2}, but is presented in a language that clearly 
involves only the properties (3.1)--(3.6) of a class $\cC$. 

\proclaim{Theorem 5.12}
Let $M$ denote a $\cC$-manifold, and let $X$ denote a 
closed $\cC$-hypersurface in $M$. Let $K$ be a compact subset
of $M$.  Then there is a neighbourhood $W$ of $K$ and a surjective 
mapping $\varphi : W' \to W$ of class $\cC$, such that: \newline
\noindent (1) $\varphi$ is a composite of finitely many 
$\cC$-mappings, each of which is either a blowing-up with 
smooth centre (that is nowhere dense in the smooth points of the
strict transform of $X$) or a surjection 
of the form 
$$
\coprod_j U_j \ \to \ \bigcup_j U_j \ , 
$$
where the latter is a finite covering of the target space by 
coordinate charts and $\coprod$ means disjoint union.\newline
\noindent (2) The final strict transform 
$X'$ of $X$ is smooth, and $\varphi^{-1} (X)$ has only 
normal crossings. 
(In fact $\varphi^{-1} (X)$ and $\det d \varphi$ simultaneously 
have only normal crossings, where $d\varphi$ is the Jacobian 
matrix of $\varphi$ with respect to any local coordinate system.)
\endproclaim 

We note two immediate consequences of Theorem 5.12. 
Let $M$ denote a $\cC$-manifold, and let $X$ denote 
any closed $\cC$-subset of $M$. 

\proclaim{Corollary 5.13}
(Rectilinearization theorem.) 
Suppose that $M$ is of (pure) dimension $n$. 
Let $K$ be a compact subset of $M$. Then there are finitely 
many mappings of class $\cC$, $\varphi_i : U_i \to M$, 
where each $U_i$ is an open neighbourhood of the origin in 
$\IR^n$, such that: \newline
\noindent (1) There is a compact subset $L_i$ of $U_i$, for 
each $i$, such that $\bigcup \varphi_i (L_i)$ is a neighbourhood
of $K$ in $M$. \newline
\noindent (2) For each $i$, $\varphi^{-1}_i (X)$ is a union of 
coordinate subspaces. 
\endproclaim 

We can obtain Corollary 5.13 by applying Theorem 5.12
locally to the hypersurface defined by the product of 
local defining equations of $X$. 

\proclaim{Corollary 5.14}
(Uniformization theorem.)  
There is a manifold $N$ of class $\cC$ and a proper 
$\cC$-mapping $\varphi : N \to M$ such that 
$\varphi (N) = X$. 
\endproclaim 

\remark{Remark 5.15}
The {\it dimension} of a closed $\cC$-subset $X$ of $M$
is well-defined (for example, by Corollary 5.14 and invariance
of domain).  If $X$ is a hypersurface, then Theorem 5.12
implies Corollary 5.14 with $\dim N = \dim X$. In general,
Corollary 5.14 follows from Corollary 5.13, 
but without the equality of dimensions. 
(Theorem 5.12 does imply that, if $X$ is a proper $\cC$-subset 
of $M$, then each $\varphi_i^{-1} (X)$ is a union of proper 
coordinate subspaces of $\IR^n$, in Corollary 5.13, and 
$\dim N < \dim M$, in Corollary 5.14.)
In order to deduce Corollary 5.14 with equality of dimensions, 
in general, from Theorem 5.12 (by the argument of \cite{BM1, Proof
of Theorem 5.1}, for example), we would need a positive answer 
to the following question: 

{\it Let $X$, $Y\subset M$ denote closed $\cC$-subsets such that 
$\dim (X\bs Y) = k$.  Is there a closed $\cC$-set $Z$ such that
$X\bs Y \subset Z$ and $\dim Z = k$ (for example, 
when $Y$ is a smooth hypersurface)?}
\endremark

\head 6. Applications\endhead

In this section, we note three applications of resolution
of singularities (or, more precisely, of the weaker version,
Theorem 5.12, and its Corollaries). 
These results seem to be new for Denjoy-Carleman classes.
Let $\cC$ denote a class of $\cC^\infty$ functions 
satisfying the hypotheses (3.1)--(3.6). 

\subhead Topological Noetherianity \endsubhead
Let $M$ be a manifold of class $\cC$, and let 
$\cO^\cC_a$ denote the local ring of germs of functions of 
class $\cC $ at a point $a\in M$. The completion of 
$\cO^\cC_a$ can be identified with the ring of formal 
power series over $\IR$ in $n$ indeterminates, where 
$n= \dim_a M$.  The following are important questions for 
Denjoy-Carleman classes: 

{\it Is $\cO^\cC_a$ Noetherian? Or, equivalently, 
is the formal power series ring flat over $\cO^\cC_a$?
(Or, is every finitely generated ideal closed in $\cO^\cC_a$?)}

The following theorem is the topological version of 
Noetherianity. 

\proclaim{Theorem 6.1}
Any decreasing sequence of closed $\cC$-subsets of $M$, 
$X_1 \supseteq X_2 \supseteq \cdots $, stabilizes
in some neighbourhood of a compact set $K$.  
(In other words, there exists $k$ such that, 
in a neighbourhood of $K$, $X_j = X_k$ for all 
$j\ge k$.)
\endproclaim 

\demo{Proof}
We can assume that $X_1 \ne M$. 
By Corollary 5.14 (and Remark 5.15), there is a proper 
$\cC$-mapping $\varphi : M_1 \to M$ such that 
$\dim M_1 < \dim M $ and $\varphi (M_1) = X_1$. Then 
$\varphi^{-1} (X_2) \supseteq \varphi^{-1} (X_3) \supseteq \cdots$
is a decreasing sequence of closed $\cC$-subsets of $M_1$, so that 
the result follows by induction on the dimension of the 
ambient manifold.
\endprf
\enddemo

\subhead \L ojasiewicz inequalities \endsubhead
Proofs of \L ojasiewicz's inequalities (Theorem 6.2 below) 
depending only on Theorem 5.12 were given in 
\cite{BM2, Section 2}. For Denjoy-Carleman classes, 
only more restrictive results, in dimension $2$, 
were previously known \cite{V}. 

\proclaim{Theorem 6.2}
I. Let $M$ denote a manifold of class $\cC$, and let $f$, 
$g\in \cC(M)$.  Suppose that $\{ x: g(x) = 0 \} \subseteq
\{ x: f(x) = 0 \}$ in a neighbourhood of a compact set $K$. 
Then there exist $c$, $\lambda > 0$ such that 
$$
| g(x) | \ \ge \ c| f(x) |^\lambda
$$
in a neighbourhood of $K$.  (The infimum of such $\lambda$
is a positive rational number.)

II.  Let $f\in \cC (U)$, where $U$ is open in $\IR^n$. 
Suppose that $K$ is a compact subset of $U$ on which 
$\grad f(x) = 0$ only if $f(x)=0$.  Then there exist
$c>0$ and $\mu$, $0<\mu \le 1$, such that 
$$
|\grad f(x) | \ \ge \ c| f(x) |^{1-\mu} 
$$
in a neighbourhood of $K$.  (Sup $\mu$ is rational.)

III. Let $f\in \cC(U)$, where $U\subset \IR^n$ is open.
Set $Z = \{ x\in U : f(x) = 0 \} $. 
Suppose that $K \subset U$ is compact. Then there are 
$c>0$ and $\nu \ge 1$ such that 
$$ 
| f(x) | \ \ge \ cd (x,Z)^\nu
$$
in a neighbourhood of $K$.  
($d (\cdot , Z)$ is the distance to $Z$.  
$\Inf \, \nu $ is rational.)
\endproclaim 

\subhead Division \endsubhead

\proclaim{Theorem 6.3}
Let $W$ be an open subset of $\IR^n$ (or a manifold of 
class $\cC$) and let $\xi \in \cC (W)$. 
Let $f\in \cC^\infty (W)$.  Suppose $f$ is {\it formally 
divisible by} $\xi $ (i.e., for all $a\in W$, $\hatf_a$ 
is divisible by $\hatxi_a$ in the ring of formal power 
series). 
Then there exists $g\in \cC^\infty (W)$ such that 
$f = \xi \cdot g$. 
\endproclaim 

\demo{Proof}
We follow Atiyah's proof of the division theorem of 
Hormander and \break
\L ojasiewicz.  (See \cite{A}, 
\cite{H\"o1}, \cite{\L}.)
Let $\varphi : W' \to W$ be a mapping of class $\cC$
as in Theorem 5.12, such that the pull-back 
$\varphi^* (\xi) := \xi \cir \varphi$ is locally a 
monomial times an invertible factor (in suitable coordinates). 
Since $\varphi^* (f)$ is formally divisible by 
$\varphi^* (\xi)$, it follows from property (3.5) 
that there is a $\cC^\infty$ function $g'$ on $W'$ such that 
$\varphi^* (f) = \varphi^* (\xi) \cdot g'$. 

Since $f$ is formally divisible by $\xi$ and $\hatvp^*_b $ is 
injective, for all $b\in W'$ (where $\hatvp^*_b$ denotes the 
formal pull-back homomorphism from Taylor series centred at 
$a = \varphi(b)$ to Taylor series centred at b), it follows 
that $g'$ is {\it formally a composite with} $\varphi$; 
i.e., for all $a\in W$, there is a formal power series $G_a$
at $a$ such that $\hatg'_b = \hatvp^*_b (G_a)$, for all
$b\in \varphi^{-1} (a)$.  Moreover, $G_a$ is uniquely determined
since $\hatvp^*_b$ is injective. 
It is enough to show there is a $\cC^\infty$ function $g$ 
on $W$ such that $\hatg_a = G_a$, for all $a$. 

Arguing inductively over the tower of mappings of which $\varphi$ 
is composed, it suffices to prove the following assertion: 
Let $U$ denote a coordinate chart of class $\cC$ in $W$, and let
$\sigma : U' \to U$ denote a blowing-up of $U$ with centre 
a coordinate subspace.  If $\eta \in \cC^\infty(U')$ is formally
a composite with $\sigma$, then there exists $\zeta \in \cC^\infty(U)$
such that $\eta = \sigma^* (\zeta)$.  This assertion is a 
special case of Glaeser's composite function theorem \cite{G}
since $\sigma$ is a very simple rational mapping.
\endprf
\enddemo

\head 7. Proof of the desingularization Theorem 5.12\endhead

We begin with a simple but important lemma on transformation of 
differential operators by blowing up (cf. \cite{H2, Section 8, (1.1)}, 
\cite{EV, Lemma 4.5}).  Consider a blowing-up $\sigma: U'\to U$, 
where $U$ is an open neighbourhood of $0$ in $\IR^n$, with 
centre a coordinate subspace $C = \{ x_i  = 0 \ , \ i\in I\}$, 
where $I \subset \{ 1,\dots,n \}$.  We use the notation of 
Definition 3.9.  Note that, for each $i\in I$, $y_\exc = y_i$ 
generates $\cI_{\sigma^{-1} (C)} $ in the chart $U_i$, and 
$$
U' \bs \bigcup\limits_{j\ne i} U_j \ = \ \big\{ y\in U_i : 
\ y_j = 0 \ , \ j\in I \bs \{ i \} \big\} \ . 
$$
The following is an easy calculation.

\proclaim{Lemma 7.1}
Let $f$ be a germ of a function of class $\cC$ at a point 
$a\in C$.  Let $e \in \IN$. 
Suppose that $\mu_{C,a} (f) \ge e$. 
Then, for each $i\in I$: \newline
\noindent (1) If $j\notin I $, then 
$$
\frac{1}{y^{e-1}_i} 
\left( \frac{\partial f}{\partial x_j} \circ \sigma \right) \ = \ y_i 
\frac{\partial}{\partial y_j} 
\left( \frac{f\circ \sigma}{y^e_i } \right) \ . 
$$
\noindent (2) If $j \in I \bs \{ i \}$, then 
$$
\frac{1}{y^{e-1}_i} 
\left( \frac{\partial f}{\partial x_j} \circ \sigma \right) \ = 
\ \frac{\partial}{\partial y_j} 
\left( \frac{f\circ \sigma}{y^e_i } \right) \ . 
$$
\noindent (3) (If $j=i$, then) 
$$
\frac{1}{y^{e-1}_i} 
\left( \frac{\partial f}{\partial x_i} \circ \sigma \right) \ = 
\ e \frac{f\circ \sigma}{y^e_i} + y_i \frac{\partial}{\partial y_i}
\left( \frac{f\circ \sigma}{y^e_i}\right) - 
\sum_{j\in I \bs \{ i \} } y_j 
\frac{\partial}{\partial y_j}
\left( \frac{f\circ \sigma}{y^e_i}\right)  \ . 
$$
\endproclaim 

\demo{Proof of Theorem 5.12}
Our aim is to define the finite sequence of transformations 
comprising the mapping $\varphi$.  At an intermediate step, 
we have both the strict transform of $X$ and the accumulated 
exceptional hypersurfaces $H_1,\dots, H_r$. Hence we consider
this more general situation from the beginning: 

Let $M$ be a manifold of class $\cC$.  Let $X$ denote a 
closed $\cC$-hypersurface in $M$, and let $H_1,\dots, H_r$
be smooth ``exceptional'' hypersurfaces in $M$, 
{\it which we do not necessarily assume have only normal 
crossings}.  Let $a\in M$.  Suppose that $s$ exceptional 
hypersurfaces pass through $a$, say $H_1,\dots, H_s$.  There is a
local coordinate chart $U$ containing $a$, with coordinates 
$(x_1,\dots , x_n)$ in which $a=0$, such that $X|U$ is defined 
by an equation of class $\cC$, 
$$ 
g(x_1,\dots, x_n) \ = \ 0 
$$
(i.e., $g\in \cC(U)$ generates $\cI_{X,x}$, for all $x\in U$). 
Write 
$$ d(x) \ := \ \mu_{X,x} \ = \ \mu_x (g) \ , \qquad 
x\in U \ . 
$$
Set $d:= d(a)$.  We can assume that $d(x) \le d$, 
for all $x\in U$, and that no exceptional hypersurfaces
other than $H_1,\dots, H_s$ intersect $U$.  After a linear 
coordinate change, we can assume that 
$$ 
\frac{\partial^d g}{\partial x^d_n } \ \ne \ 0 
{\hbox{\quad in}} \ U 
$$
and that, for each $p= 1,\dots, s$, $H_p$ is defined in 
$U$ by an equation $\lambda_p (x) = 0$, of class $\cC$, and 
$$ 
\frac{\partial\lambda_p}{\partial x_n } \ \ne \ 0 
{\hbox{\quad in}} \ U \ . 
$$

Let $z:= \partial^{d-1} g / \partial x^{d-1}_n $. 
Then $\partial z / \partial x_n \ne 0$ in $U$, so that 
$\{ z = 0\}$ defines a submanifold $N$ of $U$, of class $\cC$.
y the implicit function theorem (property (3.6)), we cansolve 
$z=0$ locally as 
$$
x_n \ = \ \varphi (x_1,\dots,x_{n-1}) \ ,
$$
where $\varphi$ is of class $\cC$. In fact, then, we can assume
that
$$
z \ = \ u \left( x_n - \varphi (x_1,\dots,x_{n-1}) \right)
$$
in $U$, where $u$ is nonvanishing and $u$ is of class $\cC$
(by property (3.5)). On course, $(x_1,\dots,x_{n-1})$ restricts
to a coordinate system on $N$; we write $\tilx = 
(x_1,\dots,x_{n-1})$ to denote this coordinate 
system. 

After a coordinate change $x'_n = x_n -\varphi (x_1,\dots, x_{n-1})$, 
$x'_j = x_j$, $j< n$, we can assume that $\varphi=0$. 
In other words, we can assume that 
$$ N \ = \ \{ z = 0 \} \ , 
{\hbox{\quad where}} 
\ z \ = \ \frac{\partial^{d-1} g}{\partial x^{d-1}_n } \ , 
$$
and that 
$$
z \ = \ u \cdot x_n \ , 
$$
where $u$ does not vanish in $U$.  In particular: \newline
\noindent (7.2) For all $x\in U$, $\mu_x (g) \ge d $ if and only if 
$x\in N$ (i.e., $x_n = 0 $) and 
$$ 
\mu_{\tilx}
\left( \frac{\partial^q g}{\partial x^q_n } 
\Big| N \right) \ \ge \ d-q \ , \quad 
q = 0 ,\dots, d-2\ . 
$$

Consider also the exceptional hypersurfaces $H_1 ,\dots, H_s$. 
Write
$$ \matrix
&c_q(g) \ := 
\ \displaystyle{\frac{\partial^q g}{\partial x^q_n}} \Big| N \ , \hfill
&q=0 ,\dots, d-2 \ , \hfill \\
&&\\
&b (\lambda_p) \ := \ \lambda_p \big| N \hfill , 
&p=1,\dots, s \ . \hfill 
\endmatrix
\tag 7.3
$$
If $x\in U$, set 
$$ 
s(x) \ := \ \# \{ p : \ x\in H_p \ , 
\ p=1,\dots , s \}
$$
(so that $s(a) = s$). Extending (7.2), we have
$$
\aligned
\{ x \in U &: \ \big( d(x), s(x)\big) = (d,s) \}
\\  
&= \{ x\in U : \ \mu_x (g) \ge d\ , 
\ \mu_x (\lambda_p) \ge 1 \ , \ p=1,\dots, s \} \\
&= \{ \tilx\in N : \ \mu_\tilx (c_q) \ge d-q\ , 
\ q=0 \ ,\dots, d-2 \ , \\
&\phantom{= \{ \tilx\in N : \ } \ \mu_\tilx (b_p) \ge 1\ , 
\ p=1 \ ,\dots, s \} \ , 
\endaligned \tag 7.4
$$
where each $c_q = c_q (g)$, $b_p = b(\lambda_p )$. 
\enddemo

\definition{Claim}
We can assume that every $c_q$ or $b_p$ that is not identically zero
satisfies
$$
\aligned
c^{d!/(d-q)}_q 
\ &= \ (\tilx^{\Omega(q)})^{d!} c^*_q \ , 
\quad q = 0 ,\dots, d-2 \ , \\
b^{d!}_p \ &= \ (\tilx^{\tau(p)})^{d!} b^*_p \ ,
\quad p = 1,\dots, s \ , 
\endaligned 
\tag 7.5
$$
where each $\Omega(q)$, $\tau(p) \in \IQ^{n-1}$, all 
$c^*_q (\tilx)$, $b^*_p (\tilx)$ are nonvanishing, and the 
collection of multiindices $\{ \Omega (q) , \tau (p)\}$ is totally 
ordered with respect to the natural partial ordering of 
$\IN^{n-1}$. (If $\Omega , \tau \in \IN^{n-1}$, then 
$\Omega \le \tau$ means $\Omega_j \le \tau_j$, 
$j=1,\dots,n-1$. $\tilx^\Omega$ denotes the monomial 
$x^{\Omega_1}_1 \cdots x^{\Omega_{n-1}}_{n-1}$.) 
\enddefinition

When the assumptions (7.5) are satisfied, we will say 
we are in the ``monomial case''.  We will prove the claim 
below, by induction on dimension.  But first we calculate 
the effect on our local equations of blowing up with suitable 
centre, since this calculation provides both the motivation 
for making the claim, and tools that we will need to complete 
the proof of the theorem once we reduce to the monomial case. 

\subhead Effect of blowing up \endsubhead
Consider a blowing-up $\sigma: U' \to U$ with centre $a$
$\cC$-submanifold $C$ of $U$ in the {\it equimuiltiple locus of}
$a=0$ (i.e., in $\{ x\in U: d(x) = d := d(a) \} $). 
Then $C \subset N$, by (7.2), so we can assume that 
$$
C \ = \ \{ \tilx = (x_1,\dots, x_{n-1}) \in N : \ x_i = 0 \ , 
\ i\in I \} \ , 
\tag 7.6
$$
where $I \subset \{ 1,\dots, n-1 \}$.  Then $U'$ is covered by 
coordinate charts $U_i $, $i\in I$, and $U_n$, as in Definition 
3.9.

Since $N = \{ x_n = 0 \}$, the strict transform $N'$ of $N$
lies in $\bigcup_{i\in I} U_i$.  For each $i\in I$, 
$$
N' \cap U_i \ = \ \{ y = (y_1,\dots, y_n): \  y_n = 0 \} \ . 
$$
(We use the notation of Definition 3.9.)

Consider $i\in I$.  Then the strict transform $X'$ of $X$ by 
$\sigma$ is defined in $U_i$ by $g' (y_1,\dots, y_n) = 0$, 
where
$$ 
g' \ = \ y^{-d}_i  g\circ \sigma\ . 
$$
By Lemma 7.1 (2), for all $q = 0,\dots, d$, 
$$
\frac{\partial^q g'}{\partial y^q_n} \ = 
\ \frac{1}{y^{d-q}_i} 
\frac{\partial^q g}{\partial x^q_n} \circ \sigma \ . 
$$
In particular, $\mu_y (g') \le d$ if and only if 
$y \in N' = \{ y_n = 0 \} $ and 
$$
\mu_{\tily}
\left( \frac{\partial^q g'}{\partial x^q_n } 
\Big| N' \right) \ \ge \ d-q \ , \quad 
q = 0 ,\dots, d-2\ . 
$$
Moreover, writing $c_q (g') := (\partial^q g' / \partial y^q_n)
| N'$, $q=0,\dots, d-2$ (cf. (7.3)), we have
$$
c_q (g') \ = \ y^{-(d-q)}_\exc c_q (g) \circ \tilsig \ , 
$$
where $\tilsig := \sigma | N' : N' \to N $; the latter is 
the blowing-up of $N$ with centre $C$. 

On the other hand, consider $a'\in U' \bs \bigcup_{i\in I} U_i$. 
Then $d(a') < d(a)$ (where $d(a') := \mu_{a'} (g')$): 
Since $a' \not\in N'$, $d(a') = d$ only if $a' \in \sigma^{-1} (C)$. 
But in the chart $U_n$, the intersection of $\sigma^{-1} (C)$
with the complement of $\bigcup_{i\in I} U_i$ is given by 
$$ \big\{ y: \ y_i = 0 \ , 
\ i\in I \cup \{ n \} \big\} \ . $$
Since $\partial^d g' / \partial y^d_n \ne 0$, it follows, 
by Lemma 7.1 (3) (applied successively with 
$(f,e) = \big(\partial^{d-1} g / \partial x^{d-1}_n \ , 1 \big) $, 
$\big( \partial^{d-2} g / \partial x^{d-2}_n \ , 2 \big) , \dots , 
(g,d)$) that if $a'\in \sigma^{-1} (C)$, then 
$g' (a') \ne 0$. 

Now consider also the exceptional hypersurfaces 
$H_1,\dots, H_s$.  Suppose that the centre $C$ of the blowing-up
$\sigma$ lies in the ``equimultiple locus of $a$ for the pair 
$\big( d(\cdot) , s(\cdot) \big)$''; i.e., 
$C \subset \{ x\in U : \big( d(x) , s(x) \big) = (d,s) \} $. 
Define $s(x')$, $x'\in U'$, analogously to $s(x)$, using the 
strict transforms $H'_p$ of the $H_p$ by $\sigma$; i.e., 
$s(x') = \# \{ p : x' \in H'_p , p=1,\dots, s\}$. 
For each $p$, $H'_p$ is defined locally by 
$\lambda'_p (y) = 0$, where 
$$ 
\lambda'_p \ = \ y^{-1}_\exc \lambda_p \circ \sigma \ . 
$$

Consider any chart $U_i $, $i\in I$ (under the assumption 
(7.6) above).  Write $b'_p = b(\lambda'_p ) := \lambda'_p \big| N'$
(cf. (7.3)).  Then 
$$ 
b'_p \ = \ y^{-1}_\exc b_p \circ \tilsig \ , \quad 
p = 1,\dots, s \ , 
$$
and $\big( d(x') , s(x') \big) = (d,s) $, where 
$x' \in U_i $, if and only if 
$$
\aligned
\mu_{\tilx'} (c'_q ) \ &\ge \ d-q \ , 
\quad q=0 ,\dots, d-2 \ , \\
\mu_{\tilx'} (b'_p ) \ &\ge \ 1 \quad\quad , \quad p=1 ,\dots, s\ .
\endaligned 
$$

\subhead To reduce to the monomial case \endsubhead
We apply the assertion of Theorem 5.12 by induction on dimension 
to (the hypersurface defined by) the function of class $\cC$ 
on the $\cC$-manifold $N$ given by the product of all 
nonzero $c^{d!/(d-q)}_q $, all nonzero $b^{d!}_p$, and all 
nonzero differences between two functions from this list. 
The claim (7.5) is then a consequence of the following 
elementary lemma \cite{BM1, Lemma 4.7}. 

\proclaim{Lemma 7.7}
Let $y=(y_1,\dots,y_m)$. 
Let $\alpha$, $\beta$, $\gamma \in \IN^m$ and let 
$a(y)$, $b(y)$, $c(y)$ be nonvanishing germs of functions of 
class $\cC$ at the origin of $\IR^m$. If 
$$
a(y) y^\alpha - b(y) y^\beta = c(y) y^\gamma \ , 
$$
then either $\alpha \le\beta $ or $\beta \le \alpha$. 
\endproclaim 

\remark{Remark 7.8}
Suppose that $\tilsig$ is a blowing-up of $N = \{ x_n = 0\}$
with smooth centre $\tilC$.  We can assume that 
$\tilC = \{ x \in N: \  x_i = 0, \ i\in I \}$, 
where $I \subset \{ 1,\dots, n-1 \}$.  Then $\tilsig$ induces
a blowing-up $\sigma $ of $U$ with centre 
$C = \{ x\in U: x_i = 0 ,\ i\in I \}$. 
In each coordinate chart $U_i$ , $i\in I$ (as defined above), 
the pull-back of $g$ (which coincides with the strict 
transform) and the $\sigma^{-1} (H_p)$ (which coincide with 
the strict tranforms $H'_p$ of the $H_p$) will retain 
the forms described above; in particular, the analogue of (7.4) 
still holds.  (Each $H'_p$ is smooth because $C$ has only 
normal crossings with respect to each $H_p$, although $C$ 
does not necessarily simultaneously have only normal crossings
with respect to the collection $\{ H_p \}$.)
Note also that the centre $C$ of the induced blowing-up 
$\sigma$ does not lie in the equimultiple locus of $a=0$. 
In these ways, Theorem 5.12 is weaker than Theorem 5.10
-- this is the price we pay to get a much simpler proof. 
\endremark

The effect of reducing to the monomial case by the inductive 
argument above is that, in addition to (the strict 
transforms of) the ``old'' exceptional hypersurfaces 
$H_1,\dots,H_s$, we have introduced ``new'' exceptional 
hypersurfaces corresponding to the blowings-up needed in the 
reduction.  This means that, in addition to $H_1,\dots, H_s$, 
we have a collection of ``new'' exceptional hypersurfaces 
that can be assumed each to be a coordinate subspace 
$x_j = 0$, where $1\le j \le n-1$. 

\subhead The monomial case \endsubhead
We assume (7.5) (and admit the possibility of other ``new'' 
exceptional hypersurfaces, each of the form 
$x_j = 0$, $1\le j \le n-1$. We consider lexicographic 
ordering of pairs $(d,s)$.  Let $S$ denote the equimultiple 
locus of $a=0$ for the pair $\big( d(\cdot), s(\cdot) \big)$; 
i.e., $S := \{ x\in U : \big( d(x), s(x) \big) = (d,s) \} $. 
Suppose that $d > 1$.
Choose coordinates as above.  

\remark{Remark 7.9}
If all $c_q \equiv 0$, then $N = \{ x \in U: \ d(x) = d \}$ and it follows
from property (3.5) that $g(x) = v(x)z^d$ in a neighbourhood of $N$, 
where $v$ is nonvanishing.
\endremark
\smallskip

Consider the case that all $c_q$ and $b_p$ vanish identically.
Let $\sigma: \ U' \to U$ be the blowing-up with centre $C = N$. If
$X'$ denotes the strict transform of $X$ by $\sigma$, then 
$X' \cap \sigma^{-1}(C) = \emptyset$; i.e., $d(x') = 0$, for all
$x' \in \sigma^{-1}(C)$.
 
Now suppose that not all $c_q$ and $b_p$ vanish identically. 
Then, by (7.4) and (7.5), 
$$
S \ = \ \{ \tilx = (x_1,\dots, x_{n-1}) \in N : 
\ \mu_\tilx (\tilx^\Omega ) \ge 1 \} \ , $$
where $\Omega := \min \{ \Omega (q) , \tau (p) \}$. 
(The meaning of the order of a monomial with rational powers is 
clear.)  Then 
$$ S \ = \ \bigcup_I \, Z_I \ , $$
where 
$$
Z_I \ := \ \{ \tilx \in N : \ x_i = 0 \ , \quad 
i\in I \} \ , 
$$
and $I$ runs over the {\it minimal} subsets of 
$\{ 1 , \dots, n-1 \}$ such that $\sum_{j\in I} \Omega_j \ge 1$; 
i.e., $I$ runs over the subsets of $\{ 1,\dots, n-1 \}$ such 
that 
$$
0 \ \le \ \sum\limits_{j\in I }  \Omega_j - 1 \ < \ \Omega_i \ , 
\quad 
{\hbox{ for all }} i \in I \ . 
$$

Consider the blowing-up $\sigma$ of $U$ with centre 
$C = Z_I$, for any such $I$ (so that $U'$ is covered by 
coordinate charts $U_i$, $i\in I \cup \{ n\}$, as before). 
In any chart $U_i $, $i\in I$, we have 
$$
\aligned
c'_q (\tily)^{d!/(d-q)} 
\ &= \ (\tily^{\Omega(q)'})^{d!} (c^*_q \circ \tilsig ) (\tily) \ , \\
b'_p (\tily)^{d!} 
\ &= \ (\tily^{\tau(p)'})^{d!} (b^*_p \circ \tilsig ) (\tily) \ , 
\endaligned
$$
where, for each $\zeta = \Omega (q) $ or $\tau (p)$, 
$$
\tily^{\zeta'} = y^{\zeta_1}_1 \cdots y^{\sum_{j\in I} \zeta_j - 1}_i 
\cdots y^{\zeta_{n-1}}_{n-1} \ . 
$$
In particular, if $a' \in \sigma^{-1} (a) \in U_i$ and 
$\big( d(a') , s(a')\big) = (d,s) = \big( d(a) , s(a) \big) $, 
then 
$$
1\ \le \ |\Omega' | \ < \ | \Omega | $$
(where $|\Omega | := \Omega_1 + \cdots + \Omega_{n-1}$). 
(Recall that if $d(a') = d$, then necessarily $a'\in U_i$, 
for some $i\in I$.) 
In particular, $\big( d(a') , s(a')\big) < 
\big( d(a) , s(a) \big) = (d,s) $ 
(throughout $U'$) after at most $d! |\Omega|$ blowings-up 
(each with centre given by a coordinate subspace of a chart
occuring as above.  Note that, by the ``monomial'' assumption, 
in each chart, $d(a') = d$ (or $s(a') = s$) at some point $a'$ 
only if $d(0) = d$ (or $s(0)= d$, respectively). 

If (in some chart) we have $d(a') = d(a)$ but $s(a') < s(a)$, 
then we can simply continue: Some $H'_p$ does not intersect
$N'$ near $a'$.  We repeat the argument of the monomial 
case above without this exceptional divisor, using the new 
$\Omega = \min \{ \Omega (q)' , \tau (p)' \}$. Finally, then, 
after finitely many such blowings-up, we have $d(a') < d(a)$
throughout each coordinate chart over $U$. 

At each step in the process, the strict transform $H'_p$ of 
each original hypersurface $H_p$ remains smooth, and each 
new exceptional hypersurface is a coordinate subspace 
$y_j = 0$, $1\le j \le n-1$.  (Thus $N'$ and the collection of 
new exceptional hypersurfaces simultaneously have only normal 
crossings.)  Moreover, each $H'_p$ and the collection 
of new exceptional hypersurfaces simultaneously have only 
normal crossings.

If $d(a) = d=1$, then $N=X$ and we can use the argument above 
to blow up until $s(a') = 0$, if $a'\in N'$.  Since the new
exceptional divisors simultaneously have only normal crossings
with respect to $N'$, we have the conclusion of the theorem 
except perhaps at points over $U$ that are outside $N' = X'$. 

It remains therefore to consider the case that $X = \emptyset$
and we have simply $s$ smooth hypersurfaces $H_1,\dots, H_s$. 
Locally, we can choose coordinates so that $H_s$ is given by 
$x_n = 0$ and each $H_p $, $1\le p< s$, is defined as before. 
(In fact, by the implicit function condition (3.6), we can assume
that, for each $p=1,\dots, s-1$, $H_p$ is defined by an equation 
of the form $x_n + b_p (x_1,\dots, x_{n-1}) = 0$.) 
The theorem now follows essentially by repeating the argument 
above, with $N=H_s$.  (For details, we refer to 
\cite{BM1, Proof of Thm. 4.4, Case 2, p. 2727}.)

\remark{Remark 7.10}
There are variants of Theorems 5.9 and 5.12 in which we avoid blowing
up with centre along which the space is ``geometrically smooth'' (or
smooth with respect to the ``reduced'' structure). Let $M$ be a manifold
of class $\cC$ and let $X$ denote a closed $\cC$-hypersurface in $M$.
Let $a \in X$, and let $g$ denote a generator of $\cI_{X,a}$. Say 
$\mu_a (g) = d$. We say that $X$ is {\it geometrically smooth at} $a$
if 
$$
g(x) \ = \ v(x)h(x)^d \ ,
$$
where $v(x), h(x)$ are of class $\cC$ and $v(a) \neq 0$. (Otherwise we
say that $a$ is a {\it geometrically singular} point.)

In terms of local coordinates as in the proof of Theorem 5.12 above,
$X$ is geometrically smooth at $a$ if and only if
$$
\displaystyle{\frac{\partial^q g}{\partial x^q_n}} \Big| N \ = \ 0 \ ,
\quad q = 0, \dots, d-2 \ ;
$$
moreover, in this case we can take
$$
h(x) \ = \ z \ := \ \displaystyle{\frac{\partial^{d-1} g}
{\partial x^{d-1}_n}}
$$
(or $h(x) = x_n$); cf. Remark 7.9. It follows that if $a$ is geometrically
singular, then $\{ x: \ \mu_x (g) \geq d \}$ contains no geometrically
smooth point near $a$.

We can obtain the following variants of Theorems 5.9 and 5.12: In the
statement of Theorem 5.9, replace the condition 
``$C_{j+1} \subset \Sing X_j$'' in (1) by ``$C_{j+1}$ lies in the
geometrically singular locus'', and replace (2) by ``$X_k$ is geometrically
smooth''. In the statement of Theorem 5.12, replace  
``centre (that is nowhere dense in the smooth points of the
strict transform of $X$)'' in (1) by 
``centre (that is nowhere dense in the geometrically smooth points of the
strict transform of $X$)'', and replace ``The final strict transform
$X'$ of $X$ is smooth'' in (2) by ``The final strict transform
$X'$ of $X$ is geometrically smooth''.

The only change needed in the proofs is to ``stop the process sooner'';
for example, in the proof above, we simply do not blow up with centre
$C = N$ in the case that all $c_q$ and $b_p$ vanish identically 
(following Remark 7.9).

It would be interesting to show that {\it the geometrically singular
locus of $X$ is a closed $\cC$-subset}.
\endremark

\enddocument